\documentclass[11pt]{article}

\usepackage{amsmath,amssymb,enumerate, enumitem,amsthm,subfigure,stmaryrd}
\usepackage{graphics,pstricks,pst-node,color,xy}
\xyoption{all}

\usepackage{tabularx}

\makeatletter
\def\hlinewd#1{%
\noalign{\ifnum0=`}\fi\hrule \@height #1 %
\futurelet\reserved@a\@xhline}

\makeatother
\numberwithin{equation}{section}
\newtheorem{thm}{Theorem}[section]

\newtheorem{dfn}[thm]{Definition}
\newtheorem*{mythm}{Theorem}

\theoremstyle{definition}

\def\BE#1{\begin{equation}\label{#1}}
\def\EE{\end{equation}}
\def\lan{\langle}
\def\ran{\rangle}

\def\blr#1{\big\lan#1\big\ran}

\def\ov#1{\overline{#1}}

\def\wt#1{\widetilde{#1}}
\def\e_ref#1{(\ref{#1})}

\def\sf#1{\textsf{#1}}
\def\tn#1{\textnormal{#1}}

\def\lra{\longrightarrow}

\def\al{\alpha}

\def\ga{\gamma}

\def\io{\iota}

\def\la{\lambda}

\def\si{\sigma}

\def\Om{\Omega}

\def\Si{\Sigma}

\def\i{\infty}
\def\hb{\hbar}

\def\C{\mathbb C}
\def\cC{\mathcal C}
\def\ctC{\wt{\mathcal{C}}}

\def\d{\mathfrak d}

\def\D{\mathfrak D}
\def\cD{\mathcal D}
\def\E{\mathbf e}

\def\nH{\textnormal{H}}

\def\cL{\mathcal L}

\def\M{\mathfrak M}

\def\cO{\mathcal O}
\def\P{\mathbb P}
\def\cP{\mathcal P}

\def\Pn{\mathbb P^{n-1}}

\def\Q{\mathbb Q}

\def\cS{\mathcal S}
\def\T{\mathbb T}

\def\V{\mathcal V}
\def\cY{\mathcal Y}
\def\Z{\mathbb Z}
\def\cZ{\mathcal Z}
\def\a{\mathbf a}

\def\cP{\mathcal{P}}
\def\bM{\mathbf{M}}
\def\nd{\textnormal{d}}

\def\ne{\textnormal{e}}

\def\ev{\textnormal{ev}}

\def\Rs#1{\underset{#1}{\mathfrak R}}

\begin{document}

\title{Localization computation of one-point disk invariants of\\ projective Calabi-Yau complete intersections}
\author{Alexandra Popa}
\date{\today}
\maketitle

\begin{abstract}
We define one-point disk invariants of a smooth projective Calabi-Yau (CY) complete intersection (CI) in the presence of an anti-holomorphic involution via localization.
We show that these invariants are rational numbers and obtain a formula for them which confirms, in particular, a conjecture by
Jinzenji-Shimizu \cite[Conjecture~1]{JS}.
\end{abstract}

\section{The one-point disk mirror theorem}
\label{intro_sec}

The problem of defining and computing disk invariants with respect to a Lagrangian
was and is the subject of much research in both mathematics and theoretical physics.
The problem of defining them is solved under certain constraints in 
\cite{So}, \cite{Cho}, and \cite{Ge1}.
The disk invariants of projective
CY CI threefolds with respect to a Lagrangian defined by the fixed locus of an anti-holomorphic involution were computed in \cite{PSW} and \cite{Sh}. 

In this paper, we define one-point disk invariants of a smooth projective CY CI in the presence of an anti-holomorphic involution via localization as
suggested by \cite[(4.7),(4.8)]{W} and \cite{JS}. Thus, our invariants are sums over graphs of rational functions in torus weights; 
see Section~\ref{dfn_subs}. We prove that these invariants are actually rational numbers and obtain a formula for them.
In the case when the CI is a hypersurface, our theorem confirms a conjecture by Jinzenji-Shimizu; see \cite[Conjecture~1]{JS}
and the \textbf{Theorem} below. 

Throughout this paper, we consider the anti-holomorphic involution
\begin{equation*}
\Om:\Pn\!\lra\!\Pn\,, \qquad
\Om([z_1,z_2,\ldots,z_n])\equiv
\begin{cases}
[\bar{z}_2,\bar{z}_1,\ldots,\bar{z}_{n},\bar{z}_{n-1}]&\hbox{if}~2|n;\\
[\bar{z}_2,\bar{z}_1,\ldots,\bar{z}_{n-1},\bar{z}_{n-2},\bar{z}_n]
&\hbox{if}~2\!\not|n.
\end{cases}
\end{equation*}
Let
$\a=(a_1,a_2,\ldots,a_l)$ be a tuple of 
positive odd integers and $$X_{\a}\subset\Pn$$
a smooth CI of multi-degree~$\a$ preserved by $\Om$ (such as the Fermat quintic threefold).
We assume that $X_{\a}$ is CY, i.e. $$\sum_{k=1}^l a_k=n$$
and that $X_{\a}$ has odd dimension, i.e.
$$n-l\in\Z^{>0}\quad\textnormal{is even}.$$

For each positive odd integer $d$, we denote by $N^{disk}_{1,d}$ the degree $d$ one-point disk invariant of $(X_{\a},\Om)$ described by Definition~\ref{disk_dfn} below.  
The invariant $N^{disk}_{1,d}$ will be expressed in terms of explicit formal power series in the variables $q^{\frac{1}{2}}$ and $q$
respectively, which we now define.
Let
\begin{equation*}\label{taudfn_e}
\tau_{\a}(q)\equiv2\sum_{\begin{subarray}{c}d\in\Z^{>0}\\d\,\textnormal{odd}\end{subarray}}q^{\frac{d}{2}}\frac{\prod\limits_{k=1}^l(a_kd)!!}{(d!!)^n}\in\Q[[q^{\frac{1}{2}}]].
\end{equation*}
The remaining formal power series necessary for our formula occur in the closed genus~$0$ and $1$ mirror formulas in \cite{bcov0}, \cite{bcov0_ci}, \cite{bcov1}, \cite{bcov1_ci}.
These are encoded by 
\begin{equation*}
F(w,q)\equiv\sum_{d=0}^{\i}
q^d\frac{\prod\limits_{k=1}^l\prod\limits_{r=1}^{a_kd}(a_kw+r)}
{\prod\limits_{r=1}^d(w+r)^n}\,.
\end{equation*}
This is a formal power series in $q$ with constant term 1 
whose coefficients are rational functions in $w$ which are regular at~$w=0$.
As in \cite{ZaZ}, we denote the subgroup of all such power series by~$\cP$
and define
\begin{equation*}
\bM:\cP\lra\cP  \qquad\hbox{by}\qquad
\bM H(w,q)\equiv \left\{1+\frac{q}{w}\frac{\nd}{\nd q}\right\}\bigg(\frac{H(w,q)}{H(0,q)}\bigg)\,.
\end{equation*}
If $p\in\Z^{\ge0}$, let
\begin{equation*}
 I_p(q)\equiv \bM^pF(0,q)\in 1+q\cdot\Q[[q]].
\end{equation*}

Let 
\BE{mirmap_e}
J(q)\equiv\frac{1}{I_0(q)}\left\{
\sum_{d=1}^{\i}q^d
\frac{\prod\limits_{k=1}^l(a_kd)!}{(d!)^n}
\left(\sum_{k=1}^{l}\sum_{r=d+1}^{a_kd}\frac{a_k}{r}\right)\right\}\in q\cdot\Q[[q]]\quad
\hbox{and}\quad
Q\equiv q\,\ne^{J(q)}.\EE
The map $q\!\lra\!Q$ is a change of variables called the \sf{mirror map}.

As is usually done in Gromov-Witten theory, we package the one-point disk invariants $N^{disk}_{1,d}$ of Definition~\ref{disk_dfn} below into a generating function
in the formal variable $Q^{\frac{1}{2}}$:
\begin{equation}\label{diskpot_e}
Z^{disk}_1(Q)\equiv\sum_{\begin{subarray}{c}d\in\Z^{>0}\\d\,\textnormal{odd}\end{subarray}}
Q^{\frac{d}{2}}N^{disk}_{1,d}.
\end{equation}

\begin{mythm}\label{disk_thm}
The one-point disk invariants $N^{disk}_{1,d}$ of $(X_{\a},\Om)$ as described by Definition~\ref{disk_dfn} below are rational numbers. They are explicitly computed by\footnotemark 
\begin{equation}\label{disk_e}
Z^{disk}_{1}(Q)=\frac{2^{\frac{n-l-2}{2}}}{I_{\frac{n-l-2}{2}}(q)}q\frac{\nd}{\nd q}\left\{\frac{1}{I_{\frac{n-l-2}{2}-1}(q)}
q\frac{\nd}{\nd q}\left\{
\frac{1}{I_{\frac{n-l-2}{2}-2}(q)}q\frac{\nd}{\nd q}\left\{\ldots q\frac{\nd}{\nd q}\left\{\frac{\tau_{\a}(q)}{I_0(q)}
\right\}\right\}\right\}\right\},
\end{equation}
where $Q$ and $q$ are related by the mirror map $Q=q\ne^{J(q)}$ in \e_ref{mirmap_e}.
\end{mythm}
\footnotetext{If $n-l-2\!=\!0$, then \e_ref{disk_e} should be read
$Z^{disk}_1(Q)\!=\!\frac{\tau_{\a}(q)}{I_0(q)}$.}
The above \textbf{Theorem} implies \cite[Conjecture~1]{JS} in the case when $X_{\a}$ is a hypersurface\footnote{See Appendix~\ref{app} for the correspondence
between the relevant notation in \cite{JS} and our notation.}.
In the case of the quintic threefold $X_{(5)}\!\subset\!\P^4$, the mirror formula \e_ref{disk_e}
recovers \cite[Theorem~1]{PSW} as can be seen using the divisor equation and $Q\frac{\nd}{\nd Q}\!=\!\frac{q}{I_1(q)}\frac{\nd}{\nd q}$.

Gromov-Witten (GW) invariants of a smooth projective variety $X$ are typically defined by integrating certain cohomology classes against a virtual fundamental class of 
a moduli space of stable maps. See \cite[Section~26.2]{MirSym} for the case of the closed invariants and \cite{Ge1} for the case of the disk invariants.
If $X$ is a CI in $\P^{n-1}$, then genus~0 Gromov-Witten invariants of $X$ are related to those of $\P^{n-1}$ via an Euler class formula
which allows computation via equivariant localization.
See \cite[Theorem~26.1.1]{MirSym} for the case of the closed invariants and \cite[Theorem~3]{PSW} for the case of the disk invariants of $X_{(5)}\!
\subset\!\P^4$. The following questions regarding the invariants $N^{disk}_{1,d}$ in Definition~\ref{disk_dfn} below remain open - in the general case
- as far as the author is aware:
\begin{itemize}
\item When and how is it possible to define these invariants as an integral over the virtual fundamental class of a moduli space of stable maps?
\item Is there an Euler class formula as suggested by \e_ref{Euler_e}, \cite[Corollary~1.10, Remark~1.11]{Ge2},  generalizing \cite[Theorem~3]{PSW}
to higher dimensions?
\end{itemize}

\subsection{Acknowledgments} I thank Aleksey Zinger for bringing \cite{JS} to my attention and for collaborating on \cite{bcov0_ci}.
The present paper is an application of one of the mirror formulas in \cite[Theorem~6]{bcov0_ci} and an extension of \cite[Lemma~6.2]{bcov0_ci}.

I thank Johannes Walcher for independently pointing out \cite{JS}.

\section{The one-point disk invariants}

\subsection{On what the invariants should be}
\label{heuristics_subs}
This section explains non-rigorously what the invariants should be  providing the motivation for Section~\ref{dfn_subs} where we actually define them. 
The idea appears in \cite{JS} and in the rich research work on disk invariants preceding it.
The notation used in this section will not be used in the remaining sections, unless re-defined.

For each positive odd integer $d$, we denote by $\ov{\M}_{0,1}(X_{\a},d)^{\Om}$ the moduli space of stable one-point degree $d$ doubled disk maps to $X_{\a}$.
An element in this space is an equivalence class $[(\Si,c),f,z_1,c(z_1)]$ where
\begin{itemize}
\item $(\Si,c)$ is a genus $0$ nodal curve together with an anti-holomorphic involution
$c\!:\!\Si\!\lra\!\Si$ so that $$\Si^{c}\equiv\left\{p\in\Si:c(p)=p\right\}$$
is a chain of circles;
\item $\Om\circ f=f\circ c$;
\item $z_1\!\in\!\Si-\Si^c$;
\item The tuple $(\Si,f,z_1,c(z_1))$ is a degree~$d$ two-point genus $0$ stable map and so it determines an element in 
$\ov\M_{0,2}(X_{\a},d)$. 
\end{itemize}

If $\Si^c$ consists of only one circle, then $\Si/c$ is a genus $0$ nodal curve with a disk attached to it.

Since $X_{\a}$ is CY, the virtual real dimension of $\ov\M_{0,1}(X_{\a},d)^{\Om}$ should be $n-l-2$.

There is a natural evaluation map $$\ev_1:\ov\M_{0,1}(X_{\a},d)^{\Om}\lra X_{\a},\qquad
\ev_1[(\Si,c),f,z_1,c(z_1)]\equiv f(z_1).$$
With $\nH\in H^2(\Pn;\Z)$ denoting the hyperplane class on $\Pn$, the one-point disk invariant that we will compute should be
of the form
\begin{equation}\label{heuristic1_e}
2\int_{\left[\ov\M_{0,1}(X_{\a},d)^{\Om}\right]^{vir}}\ev_1^*\nH^{\frac{n-l-2}{2}}.
\end{equation}
The number $2$ in \e_ref{heuristic1_e} comes from the fact that a stable one-point degree $d$ doubled disk map corresponds to two one-point disk maps,
by restricting the doubled map to either the ``upper or the lower half'' of its domain.

As in the case of closed genus~$0$ GW invariants~\cite[Section 2.1.2]{BDPP} and in that of closed reduced genus~$1$ GW invariants
\cite[Theorem~1.1]{LiZ}, we hope that
the ``one-point disk invariants'' of $(X_{\a},\Om)$ are related to those of the ambient space $\Pn$ via an Euler class formula. This is indeed the case for disk invariants without marked points of the quintic threefold; see \cite[Theorem~3]{PSW}.

We denote by $e(\V^{\Om}_{d})$ the ``Euler class of the bundle'' $$\V^{\Om}_{d}\equiv\ov\M_{0,1}^{\Om}\left
(\cL,d\right)\lra\ov\M_{0,1}^{\Om}(\Pn,d),$$
where

\begin{equation}\label{cLdfn_e}
\cL\equiv\bigoplus_{k=1}^l\cO_{\Pn}(a_k)\lra\Pn\,.
\end{equation}

Thus, we expect that \e_ref{heuristic1_e} should equal an expression of the form
\begin{equation}
\label{Euler_e}
2\int_{\left[\ov\M_{0,1}(\Pn,d)^{\Om}\right]^{vir}}e(\V^{\Om}_{d})\ev_1^*\nH^{\frac{n-l-2}{2}}.
\end{equation}
 
\subsection{Setup and closed genus~0 generating functions}

In this section, we briefly recall the setup in \cite[Sections~1.1,1.2,3]{bcov0_ci} that we need
for our graph-sum definition.

In the remaining part of this paper, all cohomology groups are with rational coefficients.

We write  $\ov\M_{0,2}(\Pn,d)$ for 
the moduli space of stable degree~$d$ maps into $\Pn$
from genus $0$ curves with $2$ marked points and
$$\ev_i:\ov\M_{0,2}(\Pn,d)\lra \Pn$$
for the evaluation map at the $i$-th marked point.
Denote by 
$$\V_{\a}\lra\ov\M_{0,2}(\Pn,d)$$
the vector bundle corresponding to the locally free sheaf
$$\pi_*\ev^*\cL\lra \ov\M_{0,2}(\Pn,d),$$
where $\cL\!\lra\Pn$ is the vector bundle giving $X_{\a}$; see \e_ref{cLdfn_e}.
For each $i\!=\!1,2$, there is a well-defined bundle map
$$\wt\ev_i\!:\V_{\a}\lra\ev_i^*\cL, \qquad
\wt\ev_i\big([\cC,f;\xi]\big)=\big[\xi(z_i(\cC))\big],$$
where $z_i(\cC)$ is the $i$-th marked point of $\cC$.
Let
$$\V_{\a}'\equiv\ker\wt\ev_1\lra \ov\M_{0,2}(\Pn,d)
\qquad\hbox{and}\qquad
\V_{\a}''\equiv\ker\wt\ev_2\lra \ov\M_{0,2}(\Pn,d).$$

We denote by $\T$ the complex $n$-torus.
The group cohomology of~$\T$ is
$$H_{\T}^*=\Q[\al_1,\al_2,\ldots,\al_n],$$
where $\al_i\!\equiv\!\pi_i^*c_1(\ga^*)$,
$\ga\!\lra\!\P^{\i}$ is the tautological line bundle,
and $\pi_i\!: (\P^{\i})^n\!\lra\!\P^{\i}$ is
the projection to the $i$-th component. We denote by
$$\Q_{\al}\equiv\Q(\al_1,\al_2,\ldots,\al_n)$$
its field of fractions.

We denote the equivariant $\Q$-cohomology of a topological space $M$
with a $\T$-action by $H_{\T}^*(M)$.
If the $\T$-action on $M$ lifts to an action on a complex vector bundle $V\!\lra\!M$,
let $\E(V)\in H_{\T}^*(M)$ denote the \sf{equivariant Euler class of} $V$.

We consider the standard action of $\T$ on $\C^{n}$,
$$\big(t_1,t_2,\ldots,t_n\big)\cdot (z_1,z_2,\ldots,z_n) 
\equiv\big(t_1z_1,t_2z_2,\ldots,t_nz_n\big).$$
This action naturally induces actions on $\Pn$ and the tautological line bundle $\ga\!\lra\!\Pn$,
on $\ov\M_{0,2}(\Pn,d)$ and $\V_{\a}$, $\V'_{\a}$, $\V''_{\a}$, and on the universal cotangent line bundles.
We denote the equivariant Euler class of 
the universal cotangent line bundle for
the $i$-th marked point by $\psi_i$.

We denote by
$$x\equiv\E(\ga^*)\in H_{\T}^2(\Pn)$$
the equivariant hyperplane class.
The equivariant cohomology of $\Pn$ is given~by
\begin{equation*}\label{pncoh_e}
H_{\T}^*(\Pn)= \Q[x,\al_1,\ldots,\al_n]\big/(x\!-\!\al_1)\ldots(x\!-\!\al_n).
\end{equation*}

With $p\!=\!\frac{n-l-2}{2}$ we define, following \cite{bcov0_ci}, 
\begin{equation}\label{Zp_e}
\cZ_p(x,\hb,Q)\equiv
x^{l+p}\!+\!\sum_{d=1}^{\i}\!Q^d\ev_{1*}\!
\left[\frac{\E(\V_{\a}'')\ev_2^*x^{l+p}}{\hb\!-\!\psi_1}\right]
\in\big(H_{\T}^*(\Pn)\big)[[\hb^{-1},Q]].
\end{equation}

\subsection{The graph-sum definition}
\label{dfn_subs}
In this section we define the one-point disk invariants $N^{disk}_{1,d}$ of $(X_{\a},\Om)$ via localization with respect to the $\T^m$-action
in \e_ref{littletorus_e} below.
Our definition is motivated by the Localization Theorem~\cite{ABo}~\cite{GraPa}, by \cite{PSW}, \cite[Lemma~3.1]{bcov0_ci}, \e_ref{Euler_e}, and
the idea of ``breaking'' a localization graph at a distinguished vertex \cite{Gi}.

For each $i\!=\!1,2,\ldots,n$, let
\begin{equation*}\label{phidfn_e}
\phi_i\equiv\prod_{k\neq i}(x\!-\!\al_k)\in H^*_{\T}(\Pn).
\end{equation*}
We denote by
$$P_1\equiv[1,0,0,\ldots,0],\qquad P_2\equiv[0,1,0,\ldots,0],\qquad P_n\equiv[0,0,\ldots,0,1]\in\Pn$$
the $\T$-fixed points in $\Pn$.
If $\eta\!\in\!H^*_{\T}(\ov\M_{0,2}(\Pn,d))$, then
\begin{equation}\label{push_e}
(\ev_1)_*\eta\big|_{P_i}=\int_{\Pn}\phi_i(\ev_1)_*\eta=\int_{\ov\M_{0,2}(\Pn,d)}\eta\ev_1^*\phi_i.
\end{equation}

Denote by $m$ the integer part of $n/2$ and by $\T^m$ the complex $m$-torus.
The embedding 
\begin{equation}\label{littletorus_e}
\io\!:\T^m\lra\T, \qquad 
(u_1,u_2,\ldots,u_m)\lra 
\begin{cases} (u_1,u_1^{-1},\ldots,u_m,u_m^{-1})&\hbox{if}~n\!=\!2m,\\
(u_1,u_1^{-1},\ldots,u_m,u_m^{-1},1)&\hbox{if}~n\!=\!2m\!+\!1,
\end{cases}
\end{equation}
induces a $\T^m$-action on $\Pn$.

Let $\la_1,\ldots,\la_m$ be the weights of the standard representation
of~$\T^m$ on~$\C^m$. It follows that
\begin{equation}\label{sp_weights_e}
(\al_1,\ldots,\al_n)\big|_{\T^m}
=\begin{cases}
(\la_1,-\la_1,\ldots,\la_m,-\la_m),&\hbox{if}~n\!=\!2m;\\
(\la_1,-\la_1,\ldots,\la_m,-\la_m,0),&\hbox{if}~n\!=\!2m\!+\!1.
\end{cases}\end{equation}

In \cite[Section~6.1]{bcov0_ci}, $\cD_{i,\ga}\!\in\!\Q(\la_1,\la_2,\ldots,\la_m)$ denotes the contribution of
the half-edge disk map whose doubled corresponds to a cover of the line through $P_i$ and $\Om(P_i)$ without a marked point to the degree $\ga$ disk invariants of the CY CI threefold. 
For each $1\!\le\!i\!\le\!2m$ and each positive odd integer $\ga$, let 
\begin{equation}\label{diskedge_e}
\cD_{1,i,\ga}\equiv\frac{\cD_{i,\ga}}{\frac{2\al_i}{\ga}}\equiv
\frac{\prod\limits_{k=1}^l(a_k\ga)!!}
{\ga\underset{(k,s)\neq(i,\ga)}{\prod\limits_{k=1}^{n}
\prod\limits_{\begin{subarray}{c}1\le s\le\ga\\ s\,\tn{odd}\end{subarray}}}
\left(s\frac{\al_i}{\ga}-\al_k\right)}
\bigg(\frac{\al_i}{\ga}\bigg)^{\frac{n\ga+l}{2}},
\end{equation}
where the above rational expression in $\al$ should be viewed in $\Q(\la_1,\la_2,\ldots,\la_m)$ via \e_ref{sp_weights_e}.

If $F\!\subset\!\ov\M_{0,2}(\Pn,r)$ is a $\T^m$-fixed component, we denote by 
$NF$ its normal bundle.

For each positive odd integer $d$, let
\begin{equation}\label{disk_dfn0}
N^{disk}_{1,d}\equiv2\sum\limits_{\begin{subarray}{c}
1\le i\le 2m\\
r\in\Z^{\ge 0}\\
\ga\in\Z^{>0},\,\ga\,\textnormal{odd}\\
2r+\ga=d\\
F\subset\ov\M_{0,2}(\P,r)^{\T^m}\,\textnormal{component}
\end{subarray}}\int_F
\frac{\E(\V''_{\a})\ev_2^*\phi_i\ev_1^*x^{\frac{n-l-2}{2}}}{(h-\psi_2)\E(NF)}
\cD_{1,i,\ga}\Big|_{\hb=\frac{2\al_i}{\ga}}\in\Q(\la_1,\la_2,\ldots,\la_m),
\end{equation}

where we set the summand equal to $\al_i^{\frac{n-l-2}{2}}\cD_{1,i,\ga}$ if $r\!=\!0$ and the sum is to be taken in $\Q(\la_1,\ldots,\la_m)$
via \e_ref{sp_weights_e}.

If $r\!>\!0$ and $F$ is a $\T^m$-fixed component in $\ov\M_{0,2}(\Pn,r)$ contributing to
\e_ref{disk_dfn0}, then the second marked point of a stable map in $F$ gets mapped to $P_i$. So on such a 
$\T^m$-fixed component,
\begin{equation*}
\E(\V''_{\a})=\frac{\E(\V'_{\a})\ev_1^*x^l}{\al_i^l}\quad\textnormal
{and so}\quad
\frac{\E(\V''_{\a})\ev_2^*\phi_i\ev_1^*x^{\frac{n-l-2}{2}}}{(h-\psi_2)\E(NF)}=
\frac{\E(\V'_{\a})\ev_1^*(x^{l+\frac{n-l-2}{2}})\ev_2^*\phi_i}{\al_i^l(\hb-\psi_2)\E(NF)}.
\end{equation*}
This together with \e_ref{Zp_e} and \e_ref{push_e} shows that \e_ref{disk_dfn0} is equivalent to the following definition.

\begin{dfn}\label{disk_dfn} With $Z^{disk}_1(Q)$ denoting the one-point disk generating function
in \e_ref{diskpot_e}, the one-point disk invariants of $(X_{\a},\Om)$
are defined by
$$Z^{disk}_1(Q)=2\sum_{1\le i\le 2m}\sum_{\begin{subarray}{c}\ga\in\Z^{>0}\\\ga\,\textnormal{odd}\end{subarray}}
Q^{\frac{\ga}{2}}\cD_{1,i,\ga}\frac{1}{\al_i^l}\cZ_{\frac{n-l-2}{2}}(\al_i,\hb,Q)\left|_{\hb=\frac{2\al_i}{\ga}}\right.,$$
where the $\T$-action is restricted to the $\T^m$-action via \e_ref{sp_weights_e}  and $\cZ_{\frac{n-l-2}{2}}$ is given by \e_ref{Zp_e}.
\end{dfn}

\section{Proof of the Theorem}
\subsection{Closed genus~$0$ mirror theorem}
In this section we recall parts of \cite[Theorem~6]{bcov0_ci} that we need for the proof of
the \textbf{Theorem.}

As in \cite{bcov0_ci}, let 
\begin{equation}\label{Ydfn_e}
\cY(x,\hb,q)\equiv\sum_{d=0}^{\i}q^d
\frac{\prod\limits_{k=1}^l\prod\limits_{r=1}^{a_kd}(a_kx+r\hb)}
{\prod\limits_{r=1}^{d}\prod\limits_{k=1}^{n}(x\!-\!\al_k\!+\!r\hb)}\in\Q(\al,x,\hb)[[q]].
\end{equation}

Following \cite[Section~3.1]{bcov0_ci}, we define
$\D^p\cY_0(x,\hb,q)$ inductively~by
\BE{Dpdfn_e}\begin{split}
\D^0\cY_0(x,\hb,q)&\equiv \cY_0(x,\hb,q)\equiv
\frac{x^l}{I_0(q)}\cY(x,\hb,q),\\
\D^p\cY_0(x,\hb,q)&\equiv
\frac{1}{I_p(q)}
\left\{x+\hb\, q\frac{\nd}{\nd q}\right\}\D^{p-1}\cY_0(x,\hb,q)\,\qquad\forall~ p\ge1.
\end{split}\EE

By \cite[Theorem~6]{bcov0_ci}, there exist $\ctC_{p,s}^{(r)}(q)\!\in\!\Q[\al_1,\ldots,\al_n][[q]]$ and $C_1(q)\!\in\!\Q[[q]]$
such that  
\begin{equation}\label{0thm_e}
\cZ_p(x,\hb,Q)=\ne^{-J(q)\frac{x}{\hb}}\ne^{-C_1(q)\frac{\si_1}{\hb}}\left\{\D^p\cY_0(x,\hb,q)
+\sum_{r=1}^{p}\sum_{s=0}^{p-r}\ctC_{p,s}^{(r)}(q)\hb^{p-r-s}\D^s\cY_0(x,\hb,q)\right\},
\end{equation}
where $p\!=\!\frac{n-l-2}{2}$, $Q$ and $q$ are related by the mirror map $Q\!=\!q\ne^{J(q)}$ in \e_ref{mirmap_e}, and $\si_1\!\equiv\!\sum\limits_{i=1}^n\al_i$.

\subsection{The proof}
In this section we prove the \textbf{Theorem} using \e_ref{0thm_e} and an extension of 
\cite[Lemma~6.2]{bcov0_ci}.

Throughout this section, we consider only the $\T^m$-action on $\Pn$ defined by \e_ref{littletorus_e} and so
in the computation below the $\T$-weights $\al_1,\ldots,\al_n$ should be expressed in terms of the $\T^m$-weights $\la_1,\ldots,\la_m$
via \e_ref{sp_weights_e}. 
We denote by
$$\Rs{z=z_0}f(z)$$
the residue of $f$ at $z_0$.

By \e_ref{sp_weights_e},
$\si_1\!=\!0$. This together with Definition~\ref{disk_dfn} and \e_ref{0thm_e} shows that
\begin{equation}\begin{split}\label{diskY_e}
Z^{disk}_1(Q)&=2\sum_{1\le i\le 2m}\sum_{\begin{subarray}{c}\ga\in\Z^{>0}\\\ga\,\textnormal{odd}\end{subarray}}
q^{\frac{\ga}{2}}\cD_{1,i,\ga}\frac{1}{\al_i^l}\Big\{\D^{\frac{n-l-2}{2}}\cY_0(\al_i,\hb,q)+\\
&\qquad\qquad\qquad\qquad
+\sum_{r=1}^{\frac{n-l-2}{2}}\sum_{s=0}^{\frac{n-l-2}{2}-r}\ctC_{\frac{n-l-2}{2},s}^{(r)}(q)\hb^{\frac{n-l-2}{2}-r-s}\D^s\cY_0(\al_i,\hb,q)\Big\}\left|_{\hb=\frac{2\al_i}{\ga}}\right.,
\end{split}\end{equation}
where $Q\!=\!q\ne^{J(q)}$.

By \e_ref{Ydfn_e} and \e_ref{Dpdfn_e},
\begin{equation}\label{Ycomp_e}
\cY_0(\al_i,\hb,q)\left|_{\hb=\frac{2\al_i}{\ga}}\right.
=\frac{\al_i^l}{I_0(q)}\sum_{d=0}^{\i}q^d 
\left(\frac{\al_i}{\ga}\right)^{nd}
\frac{\prod\limits_{k=1}^l\!\frac{(a_k(\ga+2d))!!}{(a_k\ga)!!}}
{\prod\limits_{k=1}^n
\prod\limits_{\begin{subarray}{c}\ga+2\le s\le\ga+2d\\ s~\tn{odd}\end{subarray}}
\left(s\frac{\al_i}{\ga}-\al_k\right)}.
\end{equation}

By \e_ref{diskedge_e} and \e_ref{Ycomp_e}, for each integer $p\!\ge\!0$,
\begin{equation}\begin{split}\label{sum1_e}
I_0(q)\sum_{1\le i\le 2m}\sum_{\begin{subarray}{c}\ga\in\Z^{>0}\\\ga\,\textnormal{odd}\end{subarray}}
q^{\frac{\ga}{2}}\cD_{1,i,\ga}\frac{1}{\al_i^l}\hb^p\cY_0(\al_i,\hb,q)\left|_{\hb=\frac{2\al_i}{\ga}}\right.
\qquad\qquad\qquad\qquad\\
=2^p\sum_{1\le i\le 2m}\sum_{\begin{subarray}{c}\ga\in\Z^{>0}\\\ga\,\textnormal{odd}\end{subarray}}
\sum_{\begin{subarray}{c}t\ge \ga\\t\,\textnormal{odd}\end{subarray}}
q^{\frac{t}{2}}\left(\frac{\al_i}{\ga}\right)^{p+\frac{nt+l}{2}}
\frac{\prod\limits_{k=1}^l(a_kt)!!}{
\ga\prod\limits_{k=1}^n\prod\limits_{\begin{subarray}{c}
1\le s\le t\\
s\,\textnormal{odd}\\
(k,s)\neq (i,\ga)
\end{subarray}}\left(s\frac{\al_i}{\ga}-\al_k\right)}\qquad\qquad\qquad\qquad
\\=\!
2^p\!\!\!\!\!\sum_{1\le i\le 2m}\sum_{\begin{subarray}{c}\ga\in\Z^{>0}\\\ga\,\textnormal{odd}\end{subarray}}
\sum_{\begin{subarray}{c}t\ge \ga\\t\,\textnormal{odd}\end{subarray}}\!\!
q^{\frac{t}{2}}
\Rs{z=\frac{\al_i}{\ga}}
\left\{
\frac{z^{p+\frac{nt+l}{2}}\prod\limits_{k=1}^l(a_kt)!!}{
\prod\limits_{k=1}^n\prod\limits_{\begin{subarray}{c}
1\le s\le t\\
s\,\textnormal{odd}\end{subarray}}(sz-\al_k)}
\right\}\!=\!2^p\!\!\!
\sum\limits_{\begin{subarray}{c}t\in\Z^{>0}\\t\,\textnormal{odd}\end{subarray}}\!\!
q^{\frac{t}{2}}
\Rs{w=0}
\left\{
\frac{w^{\frac{n-l-2}{2}-1-p}\prod\limits_{k=1}^l(a_kt)!!}{\prod\limits_{k=1}^n\prod\limits_{\begin{subarray}{c}
1\le s\le t\\
s\,\textnormal{odd}\end{subarray}}(s-\al_kw)}
\right\}.
\end{split}\end{equation}
The last equality above follows from the residue theorem on $\P^1$.

Note that if $\cS(x,\hb,q)\in\Q(\al,x,\hb)[[q]]$, then 
\begin{equation}\label{trick_e}
q^{\frac{\ga}{2}}\left\{\al_i+\hb q\frac{\nd}{\nd q}\right\}\cS(\al,\hb,q)\left|_{\hb=\frac{2\al_i}{\ga}}\right.=
\hb q\frac{\nd}{\nd q}\left\{q^{\frac{\ga}{2}}\cS(\al_i,\hb,q)\right\}\left|_{\hb=\frac{2\al_i}{\ga}}\right..
\end{equation}

By \e_ref{trick_e} and \e_ref{Dpdfn_e}, whenever $s\!\in\!\Z^{>0}$,
\begin{equation}\begin{split}\label{trick_e2}
I_s(q)q^{\frac{\ga}{2}}\hb^p\D^s\cY_0(\al_i,\hb,q)\left|_{\hb=\frac{2\al_i}{\ga}}\right.=
q\frac{\nd}{\nd q}\left\{q^{\frac{\ga}{2}}\hb^{p+1}\D^{s-1}\cY_0(\al_i,\hb,q)\left|_{\hb=\frac{2\al_i}{\ga}}\right.\right\}.
\end{split}\end{equation}

Using the $0\!\le\!p\!\le\!\frac{n-l-2}{2}-1$ cases of \e_ref{sum1_e}, \e_ref{trick_e2}, and induction on $s$, we obtain 
\begin{equation}\label{vanish_e}
\sum_{1\le i\le 2m}\sum_{\begin{subarray}{c}\ga\in\Z^{>0}\\\ga\,\textnormal{odd}\end{subarray}}
q^{\frac{\ga}{2}}\cD_{1,i,\ga}\frac{1}{\al_i^l}\hb^p\D^s\cY_0(\al_i,\hb,q)\left|_{\hb=\frac{2\al_i}{\ga}}\right.=0\quad\textnormal{if}
\quad p+s\le\frac{n-l-2}{2}-1.
\end{equation}

Using the $p\!=\!\frac{n-l-2}{2}$ case of \e_ref{sum1_e}, \e_ref{trick_e2}, and induction on $s$ with $0\!\le\!s\!\le\!\frac{n-l-2}{2}$, we
obtain\footnote{The $s\!=\!0$ case of the right-hand side of \e_ref{formula_e} is $2^{\frac{n-l-2}{2}}\frac{\tau_{\a}(q)}{I_0(q)}$.}
\begin{equation}\begin{split}\label{formula_e}
2\sum_{1\le i\le 2m}\sum_{\begin{subarray}{c}\ga\in\Z^{>0}\\\ga\,\textnormal{odd}\end{subarray}}
q^{\frac{\ga}{2}}\cD_{1,i,\ga}\frac{1}{\al_i^l}\hb^{\frac{n-l-2}{2}-s}\D^s\cY_0(\al_i,\hb,q)\Big|_{\hb=\frac{2\al_i}{\ga}}
\qquad\qquad\qquad\qquad\qquad\qquad\\=
\frac{2^{\frac{n-l-2}{2}}}{I_s(q)}q\frac{\nd}{\nd q}\left\{\frac{1}{I_{s-1}(q)}
q\frac{\nd}{\nd q}\left\{
\frac{1}{I_{s-2}(q)}q\frac{\nd}{\nd q}\left\{\ldots q\frac{\nd}{\nd q}\left\{\frac{\tau_{\a}(q)}{I_0(q)}
\right\}\right\}\right\}\right\}.
\end{split}\end{equation}

The \textbf{Theorem} follows from \e_ref{diskY_e}, \e_ref{vanish_e}, and the $s\!=\!\frac{n-l-2}{2}$ case of \e_ref{formula_e}.

\vspace{3mm}
\appendix
\section{Correspondence between notation in \cite{JS} and our notation}
\label{app}
While in \cite{JS}, the target space is $M_N^k\!=\!\left\{X_1^k+X_2^k+\ldots+X_N^k=0\right\}\!\subset\!\P^{N-1}$,
our target space is $X_{\a}\!\subset\!\Pn$.
\newcommand\ST{\rule[-0.5em]{0pt}{1.5em}}
$$\begin{array}{|c|c|}\hline
\cite{JS}& \textnormal{our notation}\\\hlinewd{2.3pt}
\ST \ne^x&q\\\hline
\ST \wt{L}^{k,k}_{p}(\ne^x)&I_p(q)\\\hline
\ST \ne^{t(x)}&q\ne^{J(q)}\\\hline
\ST F_0^k(x)&\textnormal{right-hand side of } \e_ref{disk_e}\\\hline
\ST h& \nH\textnormal{ in Section~\ref{heuristics_subs}}\\\hline
\ST\blr{\cO_{h^{\frac{k-3}{2}}}}_{disk,2d-1} &N^{disk}_{1,2d-1}\\\hline
\ST \tau_k(x)&\tau_{\a}(q)\\\hline
\end{array}$$

\vspace{5mm}

\noindent
{\it Department of Mathematics, Rutgers University, Piscataway, NJ 08854-8019\\
alexandra@math.rutgers.edu}\\


\begin{thebibliography}{99}

\bibitem[ABo]{ABo} M.\,Atiyah and R.\,Bott, 
{\it The moment map and equivariant cohomology}, Topology 23 (1984), 1--28.

\bibitem[BDPP]{BDPP} G.\,Bini, C.\,de Cocini, M.\,Polito, and C.\,Procesi,  
{\it On the work of Givental relative to mirror symmetry},
{\it Appunti dei Corsi Tenuti da Docenti della Scuola}, 
Scuola Normale Superiore, Pisa,~1998.

\bibitem[Cho]{Cho} C.-H., Cho,
{\it Counting real $J$-holomorphic discs and spheres in dimension four and six},
J. Korean Math. Soc.~45 (2008), no.~5, 1427--1442.

\bibitem[Ge1]{Ge1} P.~Georgieva,
{\it Orientability of moduli spaces and open Gromov-Witten invariants},
Ph.D. thesis, Stanford University, 2011.

\bibitem[Ge2]{Ge2} P.~Georgieva,
{\it The orientability problem in open Gromov-Witten theory},
math/1207.5471.

\bibitem[Gi]{Gi} A.~Givental, 
{\it Equivariant Gromov-Witten invariants}, IMRN (1996),  no.~13,
613--663.

\bibitem[GraPa]{GraPa} T.~Graber and R.~Pandharipande,
{\it Localization of virtual classes},
Invent. Math.~135 (1999), no.~2, 487--518.

\bibitem[JS]{JS} M.~Jinzenji and M.~Shimizu,
{\it Open virtual structure constants and mirror computation of open Gromov-Witten invariants of projective hypersurfaces},
math/1108.4766.

\bibitem[LiZ]{LiZ} J.~Li and A.~Zinger,
{\it On the genus-one Gromov-Witten invariants of complete intersections},
J. Differential Geom.~82 (2009), no.~3, 641--690.

\bibitem[MirSym]{MirSym} K.~Hori, S.~Katz, A.~Klemm, R.~Pandharipande, 
R.~Thomas, C.~Vafa, R.~Vakil, and E.~Zaslow, {\it Mirror Symmetry},
Clay Math.\ Inst., AMS, 2003. 
 
\bibitem[Po]{bcov1_ci} A.~Popa, 
{\it The genus one Gromov-Witten Invariants of Calabi-Yau complete intersections},
Trans. AMS 365 (2013), no. 3, 1149--1181.

\bibitem[PoZ]{bcov0_ci} A. Popa and A. Zinger,
{\it Mirror symmetry for closed, open, and unoriented Gromov-Witten invariants},
math/1010.1946.

\bibitem[PSoW]{PSW} R.~Pandharipande, J.~Solomon and J.~Walcher,  
{\it Disk enumeration on the quintic 3-fold}, J. AMS 21 (2008), 1169-1209.

\bibitem[Sh]{Sh} V.~Shende, {\it One point disc descendants of complete
intersections}, in preparation. 

\bibitem[So]{So} J.~Solomon, 
{\it Intersection theory on the moduli space of holomorphic curves 
with Lagrangian boundary conditions},
math/0606429.

\bibitem[W]{W} J.~Walcher,   
{\it Opening mirror symmetry on the quintic}, 
Comm.~Math. Phys.~276 (2007), no.~3, 671--689.

\bibitem[ZaZ]{ZaZ} D.~Zagier and A.~Zinger,  
{\it Some properties of hypergeometric series 
associated with mirror symmetry}, 
{\it Modular Forms and String Duality},  
163--177, Fields Inst.~Commun.~54, AMS, 2008. 

\bibitem[Z1]{bcov0} A.~Zinger,   
{\it Genus zero two-point hyperplane integrals in Gromov-Witten theory},
Comm.~Analysis Geom.~17 (2010), no.~5, 1--45.

\bibitem[Z2]{bcov1} A.~Zinger,
{\it The reduced genus-one Gromov-Witten invariants of Calabi-Yau hypersurfaces},  
J.~Amer.~Math.~Soc.~22 (2009),  no.~3, 691--737.


\end{thebibliography}
\end{document}